\documentclass[10pt]{amsart}

\usepackage{amssymb, amsmath, amsthm}

\usepackage{mathrsfs}
\usepackage{extarrows}
\usepackage{amsthm}
\usepackage{amssymb}

\providecommand{\n}{$\newline$}

\theoremstyle{plain}

\newtheorem*{theorem*}{Theorem}
\theoremstyle{remark}
\newtheorem*{acknowledgements}{Acknowledgements}

\newcommand{\Spec}{\operatorname{Spec}}
\newcommand{\Hom}{\operatorname{Hom}}
\newcommand{\bC}{\mathbb{C}}
\newcommand{\bF}{\mathbb{F}}
\newcommand{\bG}{\mathbb{G}}
\newcommand{\bQ}{\mathbb{Q}}
\newcommand{\bP}{\mathbb{P}}
\newcommand{\cO}{\mathcal{O}}

\newcommand{\Ext}{\operatorname{Ext}}

\usepackage{hyperref}
\hypersetup{
    colorlinks=true,       
    linkcolor=blue,          
    citecolor=blue,        
    filecolor=blue,      
    urlcolor=blue           
}

\author{Frank Gounelas}
\title{The first cohomology of separably rationally connected varieties}
\address{Frank Gounelas, Institut f\"ur Mathematik, Humboldt Universit\"at zu Berlin,
Unter den Linden 6, 10099 Berlin.}
\date{\today}

\begin{document}

\maketitle
\thispagestyle{empty}
    We call a variety rationally connected if there passes a rational curve through every
    two general points and separably rationally connected if there exists a morphism
    $f:\bP^1\to X$ such that $f^*T_X$ is an ample vector bundle. In characteristic zero
    these notions coincide, whereas they differ in characteristic $p$. Over an
    algebraically closed field of characteristic zero, a smooth projective separably
    rationally connected variety $X$ has $H^i(X, \cO_X)=0$ for $i>0$ from Hodge theory
    (see \cite[p. 249]{debarre}). In a recent preprint, Biswas and dos Santos
    \cite{bds} prove a result which easily implies that in arbitrary characteristic, at
    least $H^1(X, \cO_X)=0$. 
    \begin{theorem*}(\cite[Theorem 1.1]{bds})
        Let $X$ be a smooth projective separably rationally connected variety over $k$ an
        algebraically closed field. Let $E$ be a vector bundle over $X$ such that for each
        $k$-morphism $f:\bP^1\to X$, the pullback $f^*E$ is trivial. Then $E$ itself is
        trivial.
    \end{theorem*}
    The claim on vanishing of first cohomology can be seen as follows. Pick a class in
    $H^1(X,\cO_X) = \Ext^1(\cO_X, \cO_X)$ corresponding to a vector bundle $E$ of rank
    two. After pulling back to any $f:\bP^1\to X$, we obtain
    \begin{equation*}
        0\to \cO_{\bP^1} \to f^*E\to \cO_{\bP^1} \to 0.
    \end{equation*}
    It follows that $f^*E$ is split since $H^1(\bP^1, \cO_{\bP^1})=0$. Now from the main
    theorem in \cite{bds}, $E$ must itself be trivial. In positive characteristic we give
    another (cohomological) proof that at least $H^1(X, \cO_X)=0$, which is a special case
    of the following.

    \begin{theorem*}
        Let $X$ be a smooth projective variety over an algebraically closed field $k$ and
        $f:C\to X$ a morphism from a smooth projective curve such that $f^*T_X$ is an
        ample bundle. Then $H^1(X, \cO_X)=0$.
    \end{theorem*}

    It should be noted that not much is known about the groups $H^i(X,\cO_X)$ for smooth
    separably rationally connected varieties where $i>1$ in positive characteristic. In
    the case of smooth Fano threefolds, Shepherd-Barron \cite[Corollary $1.5$]{sb} proved
    that $H^i(X,\cO_X)=0$ for $i>0$. It is also shown (ibid.\ Corollary $12.4$) that at
    least in the case of Picard rank one, Fano threefolds are liftable to characteristic
    zero so they are separably rationally connected (in general Fano varieties are only
    rationally chain connected) and hence satisfy the conditions of the theorem above.
    Smooth separably unirational (hence separably rationally connected) threefolds have
    been shown (see \cite[Theorem $2.5$]{nygaard}) to have $H^i(X,\cO_X)=0$ for $i=1,2,3$.
    In higher dimension, Koll\'ar \cite[Theorem $11$]{kollarnonrational} has shown that
    there exist smooth Fano varieties in positive characteristic which are not even
    separably uniruled (see also \cite[Chapter V]{kollar}), yet a general Fano hypersurface
    is so by \cite[Theorem $1.4$]{zhu}. On the other hand, Fano varieties which are
    also liftable to $W_2(k)$ satisfy Kodaira vanishing by Deligne-Illusie and hence have
    $H^i(X, \cO_X)=0$ for $i>0$, but it is not known whether Fano varieties satisfy
    Kodaira vanishing (see \cite[Remark 3.5]{kollarsing}).

    \n
    For the proof of the main theorem, we proceed as follows. Similarly to the case of
    $C=\bP^1$, one proves $H^0(X, \Omega_X^m)=0$ for $m>0$ (see \cite[Proposition
    $7.4$]{gounelas}), essentially by noting that we can cover $X$ by the images of
    embeddings (see \cite[Theorem II.1.8]{kollar}) from $C$ where the restriction $T_X|_C$
    is ample. Over $\bC$ the theorem now follows as in the case of $\bP^1$ from Hodge
    theory. Note that a theorem of Bogomolov and MacQuillan (\cite{bm}, \cite{ksct}) in
    characteristic zero proves that the existence of a curve satisfying the conditions of
    the theorem implies the existence of a very free $f:\bP^1\to X$. In positive
    characteristic however this is not known, nor is it known that $X$ is rationally
    connected (see \cite{gounelas} for a discussion in this direction).  One can construct
    examples of $f:C\to X$ with $f^*T_X$ ample by starting with a very free curve
    $\bP^1\to X$ and precomposing with a finite map $C\to \bP^1$. In fact in dimension
    three and above, a general deformation of such a morphism $f$ will be an embedding
    (see \cite[Theorem II.1.8]{kollar}).
    
    \n
    The main structure of our proof in positive characteristic follows mutatis
    mutandis from the proof of Theorem $2.1$ in Nygaard's paper \cite{nygaard}. First
    note that we may restrict to the algebraic closure of a finite field by spreading out
    over a suitable finitely generated $\bF_p$-algebra $A$ and noting that $f:C\to X$
    having $f^*T_X$ as an ample bundle is an open condition in the setting over $\Spec A$.
    Hence assume $f:C\to X$ with $f^*T_X$ ample is defined over $\bar{\bF}_p$. Consider
    now the Artin-Schreier sequence of \'etale sheaves on $X$
    \begin{equation*}
        0\to\bF_p\to\bG_a\xrightarrow{F-1}\bG_a\to0.
    \end{equation*}
    The cohomologies of $\bG_a$ and $\cO_X$ agree and since the latter is coherent,
    \'etale and Zariski cohomology agree hence we may assume that all cohomology groups
    are taken in the \'etale site. We obtain an exact sequence
    \begin{equation*}
        0\to H^1(X, \bF_p)\to H^1(X,\cO_X)\xrightarrow{F-1} H^1(X, \cO_X)\to 0
    \end{equation*}
    where the last map is surjective due to SGA7.XXII Proposition $1.2$. Suppressing base
    points, we use a method of Suwa \cite{suwa} to show that a $p$-group in the \'etale
    fundamental group $\pi_1(X)$ is trivial. In the case of $C=\bP^1$ Koll\'ar has proved
    that $\pi_1(X)$ is trivial using the de Jong-Starr Theorem, see \cite[Corollaire
    3.6]{debarre} (also \cite[Remark 2.5]{bds} for a correction), although in the case of
    higher genus $C$ the \'etale fundamental group could a priori be infinite (the author
    expects this is not the case however). Suwa, using a computation in crystalline
    cohomology, first proves that the vanishing of global differential forms implies that
    $h^i_p=\dim H^i(X,\bQ_p)=0$ for $i>0$ from which $\chi_p(X) = \sum_i (-1)^ih^i_p=1$.
    The proof is identical in our setup. Note that pulling back $f:C\to X$ under an
    \'etale cover $Y\to X$ gives a smooth projective curve (possibly of higher genus)
    $g:C'\to Y$ with $g^*T_Y$ also ample. Now, let $\pi_1(X)\to G$ be any finite quotient,
    $Y\to X$ the finite \'etale cover corresponding to $G$ and let $Y\to Z$ be the degree
    $p^r$ subcover corresponding to a $p$-Sylow in $G$. From the discussion before, both
    $Y$ and $Z$ admit morphisms from curves whose pullback of the tangent bundle is ample
    and so have $\chi_p=0$. By Crew's formula, $\chi_p(Y) = p^r\chi_p(Z)$ hence $p^r =
    \deg(Y/Z)=1$. Hence $\pi_1(X)$ has no elements of order $p$ (see \cite{cl} for a
    similar argument). Now, since $\pi_1(X)$ is profinite we obtain $\pi_1^{ab}(X)
    \otimes_{\mathbb{Z}}\bF_p = 0$.
    
    \n
    Now $H^1(X, \bF_p) = \Hom(\pi_1(X), \bF_p) = \Hom(\pi_1^{ab}(X)\otimes\bF_p, \bF_p) = 0$
    and by SGA7.XXII Proposition $2.2.5$, the semi-simple component of $H^1(X, \cO_X)$
    under the endomorphism induced by Frobenius $F$ is isomorphic to $H^1(X,
    \bF_p)\otimes\bar{\bF}_p$, which is trivial. Hence $F$ is nilpotent on $H^1(X,
    \cO_X)$. The injectivity of the map of the corresponding sheaves induces $H^0(X,
    \cO_X/F\cO_X) \to H^0(X, \Omega^1_X) = 0$ and so from the cohomology of the short
    exact sequence
    \begin{equation*}
        0\to \cO_X\xrightarrow{F}\cO_X\to\cO_X/F\cO_X\to 0
    \end{equation*}
    we obtain that $F:H^1(X,\cO_X)\to H^1(X, \cO_X)$ is injective. Since $F$ is thus
    injective and nilpotent on first cohomology, the result follows.

\begin{acknowledgements}
    I would like to thank Jakob Stix for helpful conversations and to H\'el\`ene Esnault
    for some comments and bringing \cite{bds} to my attention.
\end{acknowledgements}

\bibliographystyle{alpha}

\begin{thebibliography}{Suw83}

\bibitem[BdS11]{bds} I. Biswas and J. P. P. dos Santos, {Triviality criteria for vector
bundles over rationally connected varieties}. Preprint, 2011.

\bibitem[BM01]{bm} Feodor A. Bogomolov and Michael L. MacQuillan, {Rational curves on
foliated varieties}. IHES, Preprint, 2001.

\bibitem[CL04]{cl}
Antoine Chambert-Loir.
\newblock Points rationnels et groupes fondamentaux: applications de la
  cohomologie {$p$}-adique (d'apr\`es {P}. {B}erthelot, {T}. {E}kedahl, {H}.
  {E}snault, etc.).
\newblock {\em Ast\'erisque}, (294):viii, 125--146, 2004.

\bibitem[Deb03]{debarre}
Olivier Debarre.
\newblock Vari\'et\'es rationnellement connexes (d'apr\`es {T}. {G}raber, {J}.
{H}arris, {J}. {S}tarr et {A}. {J}. de {J}ong).
\newblock {\em Ast\'erisque}, (290):Exp. No. 905, ix, 243--266, 2003.
\newblock S{\'e}minaire Bourbaki. Vol. 2001/2002.

\bibitem[Gou12]{gounelas} F. Gounelas, {\em Free curves on varieties}, preprint, arXiv:1208.4055 (2012).

\bibitem[KSCT07]{ksct}
Stefan Kebekus, Luis Sol{\'a}~Conde, and Matei Toma.
\newblock Rationally connected foliations after {B}ogomolov and {M}c{Q}uillan.
\newblock {\em J. Algebraic Geom.}, 16(1):65--81, 2007.

\bibitem[Kol95]{kollarnonrational}
J{\'a}nos Koll{\'a}r.
\newblock Nonrational hypersurfaces.
\newblock {\em J. Amer. Math. Soc.}, 8(1):241--249, 1995.

\bibitem[Kol96]{kollar}
J{\'a}nos Koll{\'a}r.
\newblock {\em Rational curves on algebraic varieties}, volume~32 of {\em
  Ergebnisse der Mathematik und ihrer Grenzgebiete. 3. Folge. A Series of
  Modern Surveys in Mathematics [Results in Mathematics and Related Areas. 3rd
  Series. A Series of Modern Surveys in Mathematics]}.
\newblock Springer-Verlag, Berlin, 1996.

\bibitem[Kol13]{kollarsing}
J{\'a}nos Koll{\'a}r.
\newblock {\em Singularities of the minimal model program}, volume 200 of {\em
  Cambridge Tracts in Mathematics}.
\newblock Cambridge University Press, Cambridge, 2013.
\newblock With a collaboration of S{\'a}ndor Kov{\'a}cs.

\bibitem[Nyg78]{nygaard}
Niels Nygaard.
\newblock On the fundamental group of a unirational {$3$}-fold.
\newblock {\em Invent. Math.}, 44(1):75--86, 1978.

\bibitem[SB97]{sb}
N.~I. Shepherd-Barron.
\newblock Fano threefolds in positive characteristic.
\newblock {\em Compositio Math.}, 105(3):237--265, 1997.

\bibitem[Suw83]{suwa}
Noriyuki Suwa.
\newblock A note on the fundamental group of a unirational variety.
\newblock {\em Proc. Japan Acad. Ser. A Math. Sci.}, 59(3):98--99, 1983.

\bibitem[Zhu11]{zhu}
Yi Zhu.
\newblock Fano hypersurfaces in positive characteristic.
\newblock Preprint, 2011.

\end{thebibliography}

\end{document}